\documentclass[preprint,review,12pt]{elsarticle}

\usepackage{amssymb}
\usepackage{amsmath}
\usepackage{amssymb}
\usepackage{amscd}
\usepackage{latexsym}
\usepackage{amsfonts}
\usepackage{graphicx}

\usepackage{euscript,amssymb,amsmath}

\newcommand{\be}[1]{\begin{equation}\label{#1}}
\newcommand{\ee}{\end{equation}}
\newcommand{\ba}[1]{\begin{eqnarray}\label{#1}}
\newcommand{\ea}{\end{eqnarray}}
\newcommand{\rf}[1]{(\ref{#1})}
\newcommand{\nn}{\nonumber}

\journal{:}

\begin{document}

\begin{frontmatter}

%% Title, authors and addresses

%% use the tnoteref command within \title for footnotes;
%% use the tnotetext command for the associated footnote;
%% use the fnref command within \author or \address for footnotes;
%% use the fntext command for the associated footnote;
%% use the corref command within \author for corresponding author footnotes;
%% use the cortext command for the associated footnote;
%% use the ead command for the email address,
%% and the form \ead[url] for the home page:
%%
%%\title{Title\tnoteref{label1}}
%% \tnotetext[label1]{}
%% \author{Name\corref{cor1}\fnref{label2}}
%% \ead{email address}
%% \ead[url]{home page}
%% \fntext[label2]{}
%% \cortext[cor1]{}
%% \address{Address\fnref{label3}}
%% \fntext[label3]{}

\title{Structural optimization of the Ziegler's pendulum: singularities and exact optimal solutions}

%% use optional labels to link authors explicitly to addresses:
%% \author[label1,label2]{<author name>}
%% \address[label1]{<address>}
%% \address[label2]{<address>}

\author{Oleg N. Kirillov }

\address{Magneto-Hydrodynamics Division (FWSH), Helmholtz-Zentrum Dresden-Rossendorf\\
P.O. Box 510119, D-01314 Dresden, Germany\\E-mail:~~o.kirillov@fzd.de}

\begin{abstract}
Structural optimization of non-conservative systems with respect to stability criteria is a research area with important applications in fluid-structure interactions, friction-induced instabilities, and civil engineering. In contrast to optimization of conservative systems where rigorously proven optimal solutions in buckling problems have been found, for non-conservative optimization problems only numerically optimized designs were reported. The proof of optimality in the non-conservative optimization problems is a mathematical challenge related to multiple eigenvalues, singularities on the stability domain, and non-convexity of the merit functional. We present a study of the optimal mass distribution in a classical Ziegler's pendulum where local and global extrema can be found explicitly. In particular, for the undamped case, the two maxima of the critical flutter load correspond to a vanishing mass either in a joint or at the free end of the pendulum; in the minimum, the ratio of the masses is equal to the ratio of the stiffness coefficients. The role of the singularities on the stability boundary in the optimization is highlighted and extension to the damped case as well as to the case of higher degrees of freedom is discussed.

\end{abstract}

\begin{keyword}
Circulatory system, structural optimization, Ziegler's pendulum, Beck's column, flutter, divergence, damping, Whitney's umbrella
%% keywords here, in the form: keyword \sep keyword

%% MSC codes here, in the form: \MSC code \sep code
%% or \MSC[2008] code \sep code (2000 is the default)

\end{keyword}

\end{frontmatter}

%%
%% Start line numbering here if you want
%%
% \linenumbers

%% main text

\section{Introduction}
\label{}

Structural optimization of conservative and non-conservative systems with respect to stability criteria is a rapidly growing research area with important applications in industry \cite{GZ88, Z89, SM03}.

Optimization of conservative elastic systems such as a problem of optimal shape of a column against buckling is already non-trivial because some optimal solutions could be multi-modal and thus correspond to a multiple semi-simple eigenvalue which creates a conical singularity of the merit functional \cite{SM03}. Despite these complications, a number of rigorous optimal solutions are known in conservative structural optimization. Nevertheless an increase in the critical divergence load given by the optimal design in such problems is usually not very large in comparison with the initial design \cite{Z89, SM03}.

In contrast to conservative systems, the non-conservative ones can loose stability both by divergence and by flutter.
It is known that mass and stiffness modification can increase the critical flutter load by hundreds percent, which is an order of magnitude higher than typical gains achieved in optimization of conservative systems \cite{C75}--\cite{KT07}. For example, Ringertz \cite{R94} reported an $838\%$ increase of the critical flutter load for the Beck's column \cite{B52} from $20.05$ for a uniform design to $188.1$ for an optimized shape. Recently Temis and Fedorov \cite{TF07} found for a free-free beam moving under the follower thrust an optimal design with the critical flutter load that exceeds that for a uniform beam by $823\%$.
We note that despite the very notion of the follower forces was debated in \cite{SL99,SL02,E05}, the Beck's column \cite{B52} as well as its discrete analogues \cite{G95, G03, L05b} including the Ziegler's pendulum \cite{Zi52} remain popular models for investigation of mode-coupling instabilities in non-conservative systems and related optimization problems.

In both conservative and non-conservative problems of structural optimization of slender structures, their optimal or optimized shapes often possess places with small or even vanishing cross-sections. The known optimized shapes of the Beck's column or of a free-free rod moving under follower thrust have almost vanishing cross-section, e.g., at the free end, which means vanishing mass of a finite element in the corresponding discretization  \cite{R94,KS98,K99,KS99,TF07,KT07}.

Another intrigue of optimization of non-conservative systems is the `wandering' critical frequency at the optimal critical load. During the optimization the eigenvalue branches experience numerous mutual overlappings and veerings \cite{C75, HW80, KK80, KS01, KS02b, KS02c} with the tendency for the critical frequency to increase and to correspond to higher modes \cite{R94, LS00a, LS00b, LS00, TF07, Bo63}. This puzzling behavior of the critical frequency still waits for its explanation.

In some problems, such as the optimal placement of the point mass along a uniform free-free rod moving under the follower thrust \cite{Sa75, PM85}, the local maxima were found to correspond to singularities of the flutter boundary such as cuspidal points where the multiple eigenvalues with the Jordan block exist \cite{K99,KS99}. In order for the last phenomenon to happen on needs at least three modes \cite{K99,KS99} meaning that the two-modes approximations \cite{Sa75, PM85} are unable to detect such optima. This reflects a general question on model reduction and validity of low-dimensional approximations in non-conservative problems discussed already by Bolotin \cite{Bo63} and Gasparini et al. \cite{G95} and recently raised up again in the context of friction-induced vibrations by Butlin and Woodhouse \cite{BW09}.

The above mentioned phenomena make rigorous proofs of the optimality in the non-conservative optimization problems substantially more difficult than in conservative ones. To the best of our knowledge, the rigorously proven optimal solutions in optimization problems for distributed circulatory systems were not found. Although in the finite-dimensional case the situation is not much better, it looks reasonable to try to understand the nature of the observed difficulties of optimization on the final dimensional non-conservative systems that depend on a finite number of control parameters.

Let us consider a circulatory system
\be{eq1}
{\bf M}\ddot {\bf x}+{\bf K}{\bf x}=0,
\ee
where dot indicates time differentiation, $\bf M$ is a real symmetric $m \times m$ mass matrix and $\bf K$ is a real non-symmetric $m \times m$ matrix of positional forces that include both potential and non-potential (circulatory) ones. The equation \rf{eq1} typically originates after linearization and discretization of stability problems for structures under follower loads, in the problems of friction-induced vibrations and even in rotor dynamics when the damping is not taken into account \cite{Bo63, KS02b, KV10}.

The characteristic equation for the circulatory system \rf{eq1} is given by the modified Leverrier's algorithm \cite{WL93}. In case when $m=2$, it reads
\be{eq2}
\det{\bf M}\lambda^4+({\rm tr}{\bf M}{\rm tr}{\bf K}-{\rm tr}({\bf M}{\bf K}))\lambda^2+\det{\bf K}=0,
\ee
where $\lambda$ is an eigenvalue that determines stability of the trivial solution.

Usually, the coefficients of the matrix $\bf K$ contain the loading parameter, say $p$, that one needs to increase by varying the coefficients of, e.g., the mass matrix. During this optimization process some masses can come close to zero so that the mass matrix can degenerate yielding $\det {\bf M}=0$. As a consequence, some eigenvalues $\lambda$ can substantially increase. On the other hand such a singular perturbation of the characteristic polynomial may cause large values of the gradient of the critical load with respect to the mass distribution. We see that a problem of optimal mass distribution in a finite-dimensional circulatory system \rf{eq1} in order to increase the critical flutter load, looks promising for explanation of the peculiarities of optimization of distributed non-conservative structures. However, it makes sense to tackle first not the most general system \rf{eq1} but rather a particular finite-dimensional non-conservative system with $m$ degrees of freedom.

In this article, we propose to take an $m$-link Ziegler's pendulum \cite{G95, G03, L05b, Zi52} as a toy model for investigation of
optimal mass and stiffness distributions that give an extremum to the critical flutter load. It appears that even the classical two-link Ziegler's pendulum had rarely been studied from this point of view in contrast to its continuous analogue --- the Beck's column. The only example of such a study known to the author is contained in the book by Gajewski and Zyczkowski \cite{GZ88}.

The paper is organized as follows. In the next Section we first consider optimization of an undamped two-link Ziegler's pendulum.
We find an explicit expression for the critical flutter load as a function of the mass or stiffness distribution and demonstrate that in the space of the two mass coefficients and the flutter load as well as in the space of the two stiffness coefficients and the flutter load, the critical flutter load forms a surface with a self-intersection and with the Whitney umbrella singular point. We consider the problem of optimal mass redistribution and find that the only two maxima of the critical flutter load correspond to a vanishing mass either in a joint or at the free end of the pendulum; in the only minimum, the ratio of the masses is equal to the ratio of the stiffness coefficients. The maxima are attained at the singular cuspidal points of the stability domain and are characterized by a dramatic increase in the critical frequency of vibrations. Then, we write down the equations of motion of an $m$-link undamped Ziegler's pendulum and consider the case $m=3$ in which we again find that the optimal mass distributions maximizing the critical flutter load correspond to vanishing of some of the three point masses. Other types of local extrema are also found that correspond to the points where the flutter boundary has a vertical tangent such as cuspidal points with triple eigenvalues, cf. \cite{K99,KS99}, or the points where the flutter boundary experiences a sharp turn, cf. \cite{KS01,KS02b,KS02c}.
Finally we consider the problem of optimal mass distribution for a two-link Ziegler pendulum with dissipation. In conclusion, we formulate an optimization problem for an $m$-link Ziegler's pendulum and discuss some hypotheses on plausible optimal solutions and their expected properties.

\section{Structural optimization of the Ziegler's pendulum}

\begin{figure}
\includegraphics[width=0.95\textwidth]{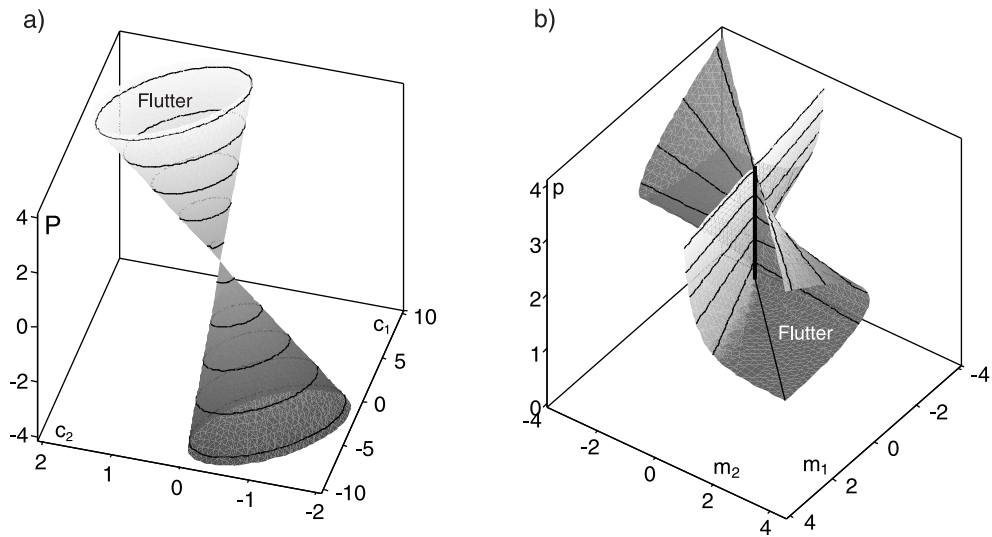}
\caption{\label{pic1} Undamped $2$-link Ziegler's pendulum. (a) The critical load $P(c_1,c_2)$ as a function of stiffness coefficients forms a conical surface in the $(c_1,c_2,P)$-space (the case when $m_1=m_2=1$ and $l=1$ is shown); (b) The critical load $p(m_1,m_2)$ as a function of masses, forms a self-intersecting surface with the Whitney umbrella singularity at the point $(0,0,2)$ of the $(m_1,m_2,p)$-space (the case when $c_1=c_2=1$ is shown). }
\end{figure}

Let us consider the classical Ziegler's pendulum consisting of two light and rigid rods of equal length $l$. The pendulum is attached to a firm basement by a viscoelastic revolute joint with the stiffness coefficient $c_1$ and the damping coefficient $d_1$. Another viscoelastic revolute joint with the stiffness coefficient $c_2$ and the damping coefficient $d_2$ connects the two rods \cite{Zi52,KV10}. At the second revolute joint and at the free end of the second rod the  point masses $m_1$ and $m_2$ are located, respectively. The second rod is subjected to a tangential follower load $P$ \cite{Zi52,KV10}.

\subsection{Undamped case}

Small deviations from the vertical equilibrium for the undamped Ziegler's pendulum are described by the equation \rf{eq1} with the mass and stiffness matrices that have the following form \cite{Zi52,KV10}
\be{z2}
{\bf M}=l^2\left(
          \begin{array}{rr}
            m_1 +m_2 & m_2  \\
            m_2  & m_2  \\
          \end{array}
        \right),\quad
{\bf K}=\left(
  \begin{array}{rr}
    c_1+c_2-Pl & Pl-c_2 \\
    -c_2 & c_2 \\
  \end{array}
\right),
\ee
where ${\bf x}=(\theta_1,\theta_2)^T$ is the vector consisting of small angle deviations from the vertical equilibrium position.

Calculating the characteristic equation $\det({\bf M}\lambda^2+{\bf K})=0$ for the Ziegler's pendulum without dissipation, we find
\be{z5}
m_1m_2l^4\lambda^4+(m_1c_2+4m_2c_2+c_1m_2-2Plm_2)l^2\lambda^2+c_1c_2=0.
\ee
By direct calculation of the roots of the characteristic equation \rf{z5} or by using the Gallina criterion \cite{G03}
we find a critical surface that separates flutter instability and marginal stability domains
\be{z6}
(2lm_2P-4m_2c_2-m_1c_1-m_2c_1)^2+(m_1c_1-m_2c_2)^2=(m_1c_1+m_2c_2)^2.
\ee
The equation \rf{z6} defines a conical surface in the $(c_1,c_2,P)$-space when $m_{1,2}$ and $l$ are fixed: flutter inside the cone, Fig.~\ref{pic1}(a).

However, in the $(m_1,m_2,P)$-space the critical surface \rf{z6} looks differently and has a form of a self-intersecting surface with the Whitney umbrella singularity at the $(0,0,2c_2/l)$-point. Indeed, expressing the critical load $P$ from \rf{z6} we get

\be{z7}
P=\frac{4m_2c_2+\left(\sqrt{m_1c_2}\pm\sqrt{m_2c_1}\right)^2}{2lm_2}\ge \frac{2c_2}{l}.
\ee
Defining the new critical flutter load as $p:=Pl/c_2$ we come to the more symmetric expression
\be{z8}
p=2+\frac{1}{2}\left(\sqrt{\frac{m_1}{m_2}}\pm\sqrt{\frac{c_1}{c_2}}\right)^2\ge 2.
\ee
The case when $c_1=c_2=1$, $m_1=2$ and $m_2=1$ corresponds to the classical result of Ziegler \cite{Zi52}
\be{z9}
p=\frac{7}{2}\pm\sqrt{2}.
\ee
The lower value of the critical load correspond to the boundary between marginal stability and flutter while the higher critical load corresponds to the conventional transition from flutter to divergence through the double zero eigenvalue with the Jordan block.

The critical load \rf{z8} as a function of the masses $p=p(m_1,m_2)$  is plotted in Fig.~\ref{pic1}(b). It is seen that the stability boundary has a self-intersection along a ray of the $p$-axis that starts at the
a Whitney umbrella singularity with the coordinates $(0,0,2)$ in the $(m_1,m_2,p)$-space. Due to the symmetry of the expression \rf{z8}, the critical load as a function of the stiffness coefficients, $p=p(c_1,c_2)$, forms the identical surface in the $(c_1,c_2,p)$-space. In the following we will study the function $p=p(m_1,m_2)$.

\begin{figure}
\includegraphics[width=0.9\textwidth]{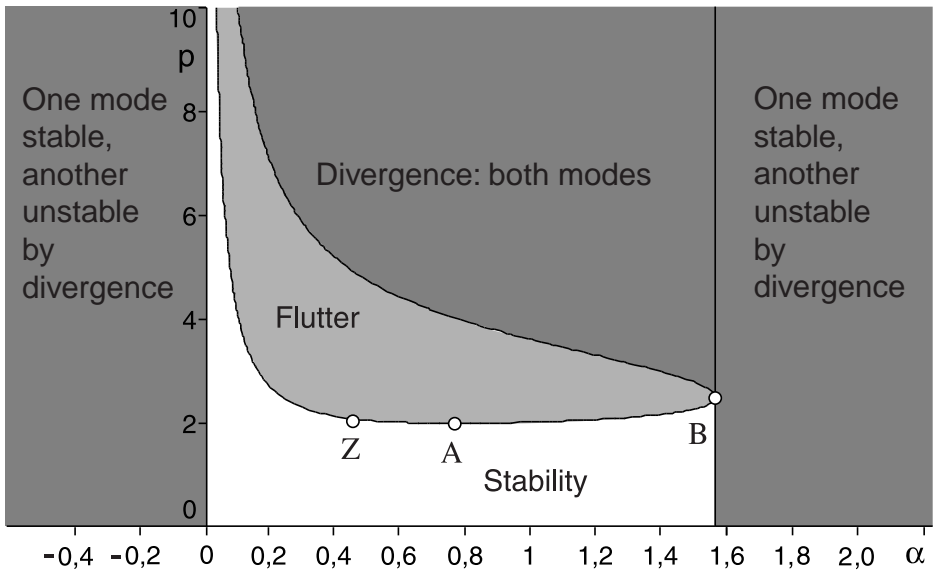}
\caption{\label{pic2} Undamped $2$-link Ziegler's pendulum with $c_1=1$, $c_2=1$: The critical flutter load $p(\alpha)$ as a function of the azimuth angle $\alpha$ indicating the direction in the $(m_1, m_2)$-plane. The point $A$ is an absolute minimum of the flutter load: $p_A=2$, the point $B$ corresponds to the local maximum: $p_B=2+c_1/(2c_2)$ with $m_1=0$, and the absolute maximum corresponds to a point $C$ (not shown) with $p_C=+\infty$ and $m_2=0$. The point $Z$ corresponds to the Ziegler's original design: $m_1/m_2=2$.}
\end{figure}

According to the inequality \rf{z8} the critical load is always not less than $p_0=2$. The minimum is reached when the masses satisfy the constraint
\be{z10}
m_1c_2=m_2c_1.
\ee
Note that the equal stiffness coefficients $c_1=c_2$ imply equal masses $m_1=m_2$. This situation corresponds to a uniformly distributed mass and stiffness in continuous systems such as the Beck's column \cite{B52}.

Usually, in the structural optimization problems the uniformly distributed stiffness and mass are considered as the initial design that is a starting point in optimization procedures. The critical load of the optimized structure is conventionally compared to that of the same structure with the uniform mass and stiffness distributions \cite{C75, Sa75, HW80, KK80, PM85, R94, KS98, LS99, LS00a, LS00b, LS00, TF07, KT07}.

Since $p(m_1,m_2)$ is a ruled surface and thus $p$ effectively depends on the mass ratio only, it is convenient to introduce the azimuth angle $\alpha$ by assuming $m_1=\cos\alpha$ and $m_2=\sin\alpha$ and to plot the critical load as a function of $\alpha$. In Fig.~\ref{pic2} the curves $p=p(\alpha)$ bound the flutter domain shown in the light gray. When $\alpha$ tends to zero, which corresponds to the vanishing mass $m_2$, the critical load increases to infinity. When $\alpha$ tends to $\frac{\pi}{2}$ and, correspondingly, the mass $m_1$ is vanishing, then the critical flutter load increases to the value
\be{z11}
p_B=2+\frac{1}{2}\frac{c_1}{c_2}.
\ee
At the point $B$ in the stability diagram of Fig.~\ref{pic2} the flutter boundary has a vertical tangent, which is a typical phenomenon in non-conservative optimization \cite{KS01, KS02b, KS02c}.

The lower part of the flutter boundary corresponds to the designs with a complex conjugate pair of pure imaginary double eigenvalues with the Jordan block; the upper part is the designs with two real double eigenvalues of the same magnitude and different sign, each with the Jordan block. Above the upper flutter boundary lies the domain of statical instability or divergence with two unstable modes corresponding to two different positive real eigenvalues. At the point $B$ the flutter boundary is tangent to the vertical part of the divergence boundary. To the right of this vertical line there is a pair of pure imaginary eigenvalues and a pair of real eigenvalues with the same magnitude and different signs. Transition from the stability boundary to the divergence boundary below the point $B$ happens when a pair of pure imaginary eigenvalues goes out of the origin in the complex plane to infinity, merge there and return back along the real axis. This happens because at the boundary $m_1=0$, i.e. $\det{\bf M}=0$. Similar divergence boundary exists at $\alpha=0$ that corresponds to $m_2=0$. Transition through the vertical line above the point $B$ is accompanied by another eigenvalue bifurcation at infinity: two real eigenvalues of the same magnitude and different signs go out of the origin in the complex plane in order to merge at infinity and then come back along the imaginary axis. The very point $B$ corresponds to an antagonist of a quadruple zero eigenvalue with the Jordan block, i.e. to a quadruplet of complex eigenvalues that merge at infinity in the complex plane.

To summarize, the popular initial design corresponding to uniformly distributed mass and stiffness turns out to give an absolute minimum to the critical flutter load of the Ziegler's pendulum. The critical flutter load attains its local maximum, $p_B$, for $m_1=0$ at the singular cuspidal point $B$ of the stability boundary where the flutter domain has a vertical tangent and touches the boundary of the divergence domain. Note that in \cite{K99} a local extremum of the flutter load for the free-free beam carrying a point mass was found to be at the cuspidal point on the flutter boundary too. The global maximum of the critical flutter load for the undamped Ziegler's pendulum is at infinity when $m_2=0$.

The global maximum corresponds to a vanishing mass at the free end of the column which qualitatively is in agreement with the numerically found optimized designs of the Beck's column available in the literature \cite{R94, LS99, LS00a, LS00b, LS00, TF07, KT07}. Indeed, all known optimized designs of the Beck's column are characterized by the vanishing cross-sections at the free end. Moreover, the gradients of the critical flutter load with respect to the mass or stiffness distribution of the Beck column are large, which is, again, in qualitative agreement with our stability diagram of Fig.~\ref{pic2}. The most interesting is the fact that with the increase of the critical flutter load the higher and higher modes were reported to be involved into the coupling that indicates the onset of flutter \cite{R94, LS99, LS00a, LS00b, LS00, TF07, KT07}. Our simple model shows that this phenomenon seems to be natural for the optimal design that causes the degeneracy in the mass matrix that gives rise to the critical frequency that increases without bounds. It would be interesting to derive the optimal solutions for the optimization of the Ziegler's pendulum by means of the Pontryagin's maximum principle, which is used in numerical optimization of the distributed non-conservative systems, see e.g. \cite{KS98, K99}.

\subsection{The $m$-link Ziegler's pendulum}

It would be very much desirable to extend our study to the case of the multiple-degrees-of-freedom Ziegler's pendulum.
The corresponding models and recursive schemes for deriving the equations of motions where proposed by Gasparini et al. \cite{G95} and by Lobas \cite{L05b}. A three-link Ziegler's pendulum was considered by Gallina \cite{G03}.

We take the linearized equations of Lobas \cite{L05b} for its model deals with the arbitrary masses and stiffnesses in the joints of an $m$-link Ziegler's pendulum in contrast to the models of Gasparini and Gallina \cite{G95,G03}.
The mass and stiffness matrices of the Ziegler-Lobas model look like
\be{l1}
{\bf M}=l^2\left(
     \begin{array}{lllll}
       \sum_{i=1}^m m_i & \sum_{i=2}^m m_i & \cdots & \sum_{i=m-1}^m m_i & \sum_{i=m}^m m_i \\
       \sum_{i=2}^m m_i & \sum_{i=2}^m m_i & \cdots & \sum_{i=m-1}^m m_i & \sum_{i=m}^m m_i \\
       \cdots & \cdots & \cdots & \cdots & \cdots\\
       \sum_{i=m-1}^m m_i & \sum_{i=m-1}^m m_i & \cdots & \sum_{i=m-1}^m m_i & \sum_{i=m}^m m_i \\
       \sum_{i=m}^m m_i & \sum_{i=m}^m m_i & \cdots & \sum_{i=m}^m m_i & \sum_{i=m}^m m_i \\
     \end{array}
   \right),
\ee
\be{l2}
{\bf K}=\left(
          \begin{array}{rrrrr}
            c_1+c_2-Pl & -c_2 & \cdots & 0 & Pl \\
            -c_2 & c_2+c_3-Pl & \cdots & 0 & Pl \\
            \cdots & \cdots & \cdots & \cdots & \cdots\\
            0 & 0 & \cdots & c_{m-1}+c_m-Pl & -c_{m-1}+Pl\\
            0 & 0 & \cdots & -c_{m-1} & c_m\\
          \end{array}
        \right).
\ee
For $m=2$ the matrices \rf{l1} and \rf{l2} are reduced to \rf{z2}. Note that
\be{l1a}
\det{\bf M}= l^{2m}\prod_{i=1}^m m_i.
\ee

\begin{figure}
\includegraphics[width=0.85\textwidth]{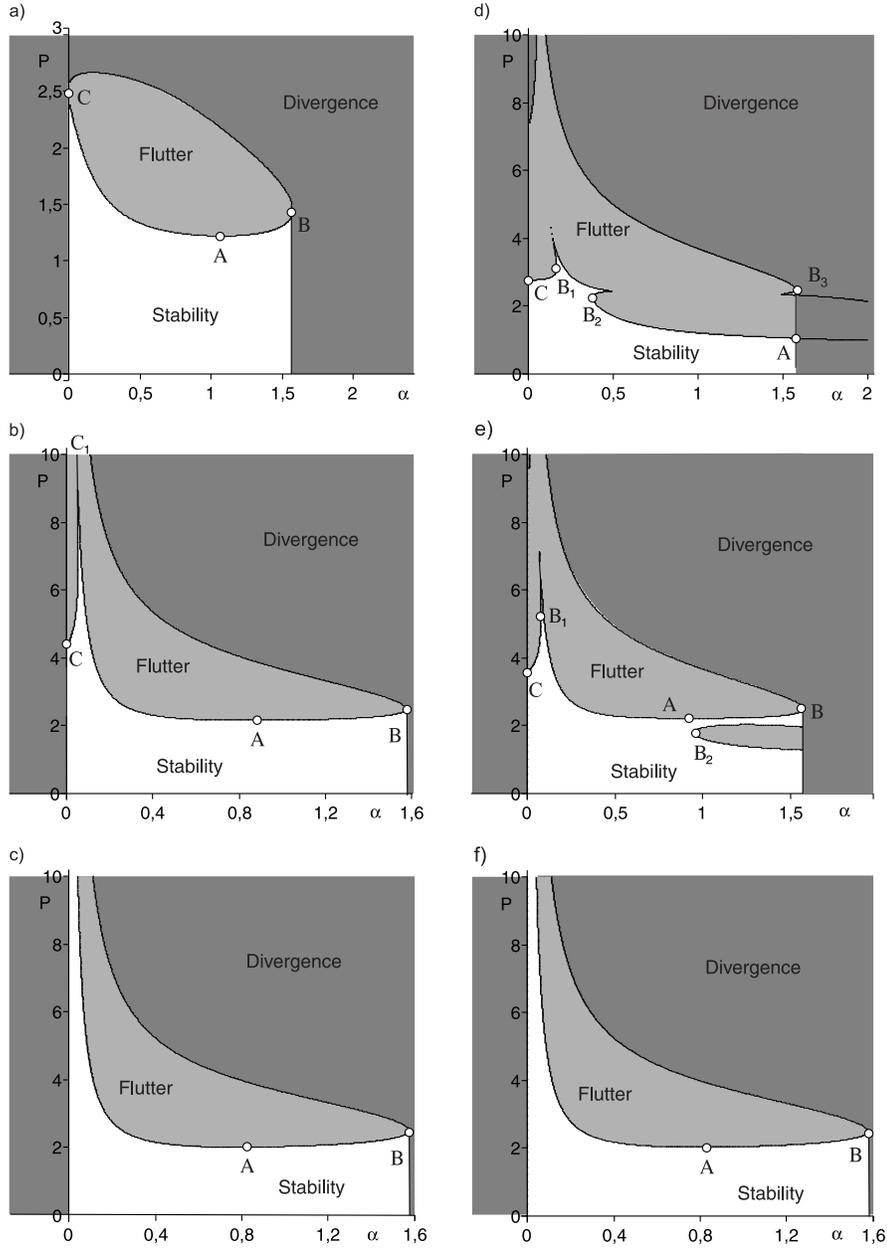}
\caption{\label{pic3} Undamped $3$-link Ziegler's pendulum with $l=1$, $c_1=c_2=c_3=1$: The critical flutter load $P(\alpha)$ as a function of the azimuth angle $\alpha$ indicating the direction in the $(m_2, m_3)$-plane at the fixed radial distance $r=1$ for (a) $m_1=0$, (b) $m_1=10$, and (c) $m_1=200$ and for the fixed $m_1=5$ and (d) $r=1$, (e) $r=0.65$ and (f) $r=0.1$.}
\end{figure}

Let us consider the Ziegler-Lobas pendulum with $m=3$ links.
Now the mass at the free end is $m_3$, while the masses $m_2$ and $m_1$ are located at the third and second joints, respectively.
The first joint connects the first rod with the basement. The length of each of the three rods is equal to $l$. The stiffness coefficients of the joints are $c_3$, $c_2$, and $c_1$, respectively.

For simplicity we assume that $c_1=c_2=c_3=c$. Then, the characteristic polynomial has the form
\be{l3}
a_0\lambda^6+a_1\lambda^4+a_2\lambda^2+a_3=0,
\ee
with the coefficients
\ba{l4}
a_0&=&l^6m_1m_2m_3,\nn\\
a_1&=&cl^4(6m_2m_3+5m_1m_3+m_1m_2)-2l^5Pm_3(m_1+m_2),\nn\\
a_2&=&3P^2l^4m_3-2(7m_3+m_2)Pl^3c+(m_1+5m_2+14m_3)l^2c^2,\nn\\
a_3&=&c^3.
\ea
Composing the discriminant matrix \cite{G03}
\be{l5}
{\bf S}=\left(
  \begin{array}{cccccc}
    a_0 & a_1 & a_2 & a_3 & 0 & 0 \\
    0 & 3a_0 & 2a_1 & a_2 & 0 & 0 \\
    0 & a_0 & a_1 & a_2 & a_3 & 0 \\
    0 & 0 & 3a_0 & 2a_1 & a_2 & 0 \\
    0 & 0 & a_0 & a_1 & a_2 & a_3 \\
    0 & 0 & 0 & 3a_0 & 2a_1 & a_2 \\
  \end{array}
\right)
\ee
and calculating the discriminant sequence consisting of the determinants of the three main minors of even order, we find that
$\Delta_1=3l^{12}m_1^2m_2^2m_3^2>0$. The expressions for determinants $\Delta_2$ and $\Delta_3=\det{\bf S}$ are rather involved and by this reason we omit them here. However, the numerical experiments evidence that the stability boundary is given by the equation $\Delta_3=0$ for the stability condition $\Delta_3>0$ implies $\Delta_2>0$.

In Fig.~\ref{pic3} using the inequality $\Delta_3>0$ we present the stability diagrams in the $(\alpha,P)$-plane where the azimuth angle $\alpha$ in the $(m_2,m_3)$-plane is introduced by assuming $m_2=r\cos\alpha$ and $m_3=r\sin\alpha$. We assume
the equal lengths of the links $l=1$ and the equal stiffness coefficients $c_1=c_2=c_3=1$ and vary the radial distance $r$ in the $(m_2,m_3)$-plane and the mass $m_1$.

Since for $m=3$ the critical surface $P(m_2,m_3)$ is no more a ruled surface as it was in the case $m=2$, the pictures in the $(\alpha,P)$-plane change with the variation of the radial distance $r$ that complicates the optimization problem.
Nevertheless such diagrams are convenient for the analysis of the geometry of the stability boundary and thus for the identification of potential extrema. Moreover, the critical surface $P(m_2,m_3)$ has a self-intersection along a ray of the
$P$-axis that starts from the singularity Whitney umbrella as in the case of the $2$-link pendulum. Therefore, at small values of $m_2$ and $m_3$ the critical load can be locally approximated by a ruled surface.

In the left column of Fig.~\ref{pic3} the radial distance $r$ in the $(m_2,m_3)$-plane is fixed to $r=1$ while the mass $m_1$ is increasing. As in the case $m=2$, (marginal) stability is possible for $\alpha \in [0,\pi/2]$. For $m_1=0$, two finite maximal values $P_A$ and $P_B$ are identified at $\alpha=0$ $(m_3=0)$ and $\alpha=\pi/2$ $(m_2=0)$, respectively, Fig.~\ref{pic3}(a).
Both maxima are attained at the cuspidal points of the stability boundary where it has the vertical tangents.
However, the stability diagram changes when $m_1=10$, Fig.~\ref{pic3}(b). Again, local extrema exist at the boundary points
$\alpha=0$ $(m_3=0)$ and $\alpha=\pi/2$ $(m_2=0)$, while the global maximum is at the point of the sharp turn of the flutter boundary with the vertical tangent near the cuspidal point $C_1\simeq(0.0403477,11.961144)$ that corresponds to triple pure imaginary eigenvalues $\lambda\simeq \pm i1.1635243$ with the Jordan block of third order, cf. \cite{K99,KS99} where at such a singular point a maximum of the critical flutter load was found for a free-free beam under the follower thrust.

With the further increase in the first mass up to $m_1=200$, the stability diagram converges to that similar to the diagram of the two-link pendulum, cf. Fig.~\ref{pic2} and Fig.~\ref{pic3}(c). This is not surprising because big inertia of the first joint makes the three-link pendulum effectively a two-link one.

The right column in Fig.~\ref{pic3} corresponds to the fixed first mass $m_1=5$ and varied radial distance $r$ in the $(m_2,m_3)$-plane. Small values of $r$ correspond effectively to a two-link pendulum. That is why Fig.~\ref{pic3}(f) with $r=0.1$ looks similar to Fig.~\ref{pic3}(c) and Fig.~\ref{pic2}.

Increase of $r$ is accompanied by a complication of the stability diagram. In particular, two cusp point singularities originate corresponding to triple pure imaginary eigenvalues with the Jordan block, Fig.~\ref{pic3}(d,e). Near these singularities, the stability boundary experience a sharp turn at the points $B_1$ and $B_2$ where the tangent to the boundary is vertical. At such points the eigencurves ${\rm Im}\lambda(P)$ form crossing that can change either to avoided crossing or to the overlapping with the origination of a bubble of complex eigenvalues in dependence on the direction of variation of the azimuth angle $\alpha$.
This phenomenon was observed in many numerical studies of non-conservative optimization problems \cite{C75,HW80,KK80,R94,KS98,LS99,LS00a,TF07,KT07} and had been described analytically by Kirillov and Seyranian in \cite{KS01,KS02b}.
Moreover, the critical load at these points can jump to a higher value corresponding to the merging of other (often higher) modes and is thus undefined \cite{C75,HW80,KK80}. Nevertheless, these points could be local extrema of the merit functional, see \cite{KS02c} where the necessary conditions for that were derived.

We stress that due to finite number of degrees of freedom and finite number of control parameters the stability boundary of an $m$-link Ziegler's pendulum can be thoroughly analyzed both analytically and numerically. In particular, the coordinates of the singular points can be calculated with arbitrary precision and thus the issues of high sensitivity of the critical flutter load at the optima could be successfully resolved in this model, contrary to the optimization problems for distributed systems. A possibility to work with the singularities related to coalescence of more than two eigenvalues allows one to investigate qualitatively the question `Should low-order models be believed?' \cite{BW09}. Indeed, two-mode approximations work well far from such singularities, in their vicinity one has to take into account higher modes.

\begin{figure}
\includegraphics[width=0.95\textwidth]{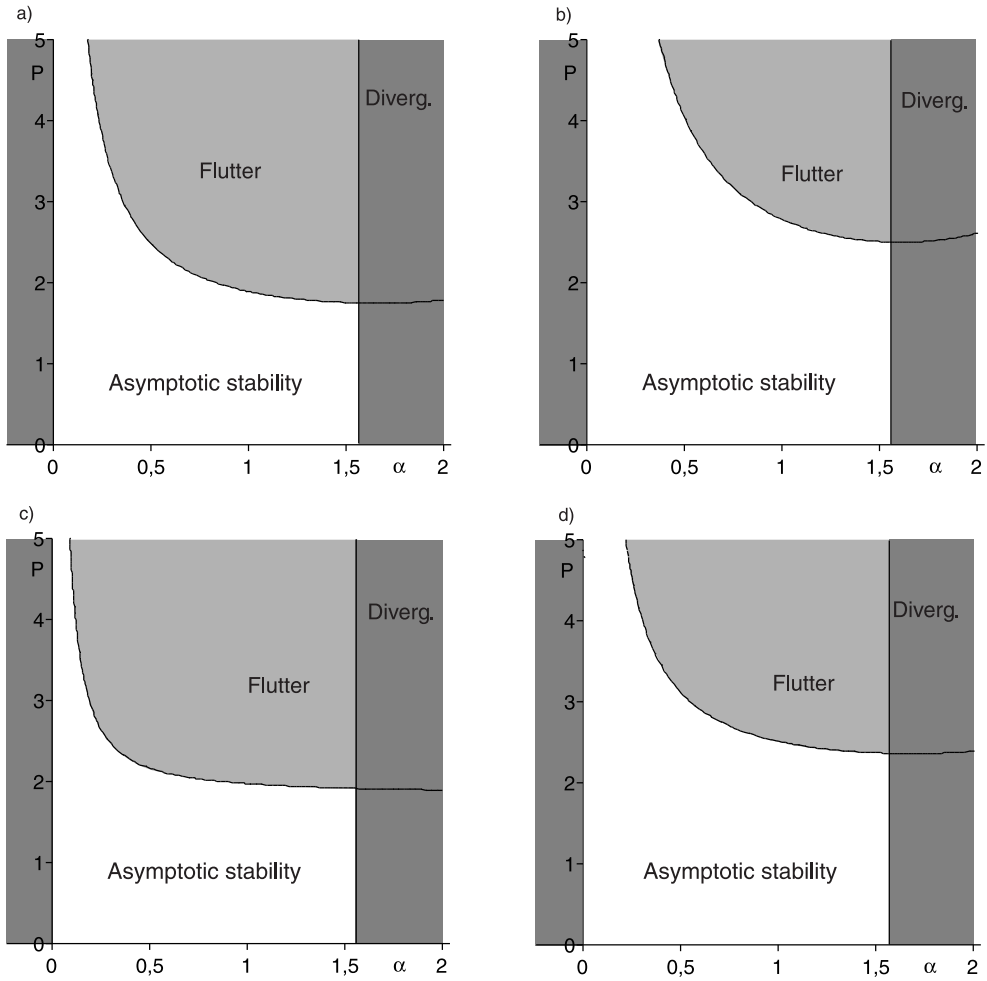}
\caption{\label{pic4} Damped $2$-link Ziegler's pendulum with $l=1$, $c_1=c_2=1$: The critical flutter load $P(\alpha)$ as a function of the azimuth angle $\alpha$ indicating the direction in the $(m_1, m_2)$-plane for (a) $d_1=d_2=1$ and $r=1$,
(b) $d_1=d_2=1$ and $r=0.4$, (c) $d_1=1$, $d_2=0.1$ and $r=1$, (d) $d_1=1$, $d_2=0.1$ and $r=0.1$.}
\end{figure}

\subsection{Damped case}

In the presence of dissipation the equation of the Ziegler's pendulum is
\be{d1}
{\bf M}\ddot {\bf x}+{\bf D}\dot{\bf x}+{\bf K}{\bf x}=0
\ee
with the damping matrix \cite{L05b, Zi52}
\be{d2}
{\bf D}=\left(
          \begin{array}{rrrrr}
            d_1+d_2 & -d_2 & \cdots & 0 & 0 \\
            -d_2 & d_2+d_3 & \cdots & 0 & 0 \\
            \cdots & \cdots & \cdots & \cdots & \cdots\\
            0 & 0 & \cdots & d_{m-1}+d_m & -d_{m-1}\\
            0 & 0 & \cdots & -d_{m-1} & d_m\\
          \end{array}
        \right).
\ee

For the two-link damped Ziegler's  pendulum with $m=2$, we find the characteristic equation
\be{d3}
a_0\lambda^4+a_1\lambda^3+a_2\lambda^2+a_3\lambda+a_4=0.
\ee
with the coefficients
\ba{d4}
a_0&=&l^4m_1m_2,\nn \\
a_1&=&l^2(m_1d_2+d_1m_2+4m_2d_2),\nn \\
a_2&=&d_1d_2+m_1l^2c_2+4m_2l^2c_2+c_1m_2l^2-2Pl^3m_2,\nn \\
a_3&=&d_1c_2+c_1d_2,\nn \\
a_4&=&c_1c_2.
\ea
Applying the Routh-Hurwitz criterion, we find the critical flutter load
\be{d5}
 P= \frac{4m_2^2(d_2^2c_1^2+d_1^2c_2^2)+d_1d_2(8m_2(m_1+2m_2)c_2^2+(c_1m_2-m_1c_2)^2)}
 {2(m_2l(4m_2d_2+d_1m_2+m_1d_2)(c_1d_2+d_1c_2)}+\frac{1}{2}\frac{d_1d_2}{m_2l^3}.
\ee
For $c_1=c_2=1$, $l=1$, $m_1=2$ and $m_2=1$ it was found to be \cite{HJ65}
\be{d6}
 P= \frac{4d_1^2+33d_1d_2+4d_2^2}
 {2(6d_2+d_1)(d_2+d_1)}+\frac{1}{2}{d_1d_2}.
\ee
The equation \rf{d6} defines a surface with the Whitney umbrella singularity in the $(d_1,d_2,P)$-space  which explains Ziegler's destabilization paradox by vanishing dissipation \cite{Zi52} as it was first demonstrated by Bottema already in 1956 \cite{KV10,B56}.

In contrast to the previous studies, e.g. \cite{KV10,HJ65,B56,Ki04}, we consider the critical flutter load \rf{d5} as function of masses $P=P(m_1,m_2)$ for the fixed damping distribution.

In Fig.~\ref{pic4}, the stability diagrams for the damped two-link Ziegler's pendulum are shown in the assumption of $c_1=c_2=1$, $l=1$, $m_1=r\cos\alpha$ and $m_2=r \sin\alpha$. The critical load $P(m_1,m_2)$ does not constitute a ruled surface and thus the function $P(\alpha)$ depends on the radial distance $r$ in the $(m_1,m_2)$-plane. We see that the extrema again correspond to the boundary points $m_1=0$ and $m_2=0$ although the singularity at $\alpha=\pi/2$ is an intersection point. Nevertheless, with the increase of the number of degrees of freedom and control parameters new types of singularities will appear. The planar stability diagrams of Fig.~\ref{pic4} depend on the damping distribution but do not tend to that of the undamped pendulum when damping goes to zero (destabilization paradox \cite{KV10,B56,Ki04}).

\subsection{`Problema novum ad cuius solutionem Mathematici invitantur'}

`A new problem that mathematicians are invited to solve' is a translated from Latin title of 1696 work by J. Bernoulli where he proposed the famous Brachistochrone problem \cite{B96}. Supporting this good old tradition we would like to formulate the following optimization problem:

{\em Given a circulatory system \rf{eq1} with the matrices $\bf M$ and $\bf K$ defined as in \rf{l1} and \rf{l2}. Find all local extrema, the absolute maximum of the critical flutter load $P$, and the corresponding extremal mass distributions $\{m_1,m_2,\ldots,m_m\}$.}

One can consider also the problem of optimal stiffness distribution or even vary both stiffness and mass. The same problems can be formulated for the damped pendulum with the damping matrix \rf{d2}.

We expect that both in the undamped and in the damped case there exists a class of extrema, corresponding to the distributions $\{m_1,m_2,\ldots,m_m\}$ with some masses $m_i=0$. It should be possible to find these optimal mass distributions explicitly perhaps with the use of the Pontryagin's maximum principle. It would be interesting to identify the singularities of the stability boundary that correspond to these extrema; some of them should be related to infinite eigenvalues $\lambda$.

On the other hand, some local extrema should exist with the mass distributions that do not contain the vanishing masses $m_i$. It would be interesting to understand at which points --- smooth or non-smooth --- of the stability boundaries they are attained. Since the system is finite-dimensional and contains the finite number of control parameters with the clear physical meaning, the locations of the singularities corresponding to multiple eigenvalues can easily be found numerically with the high accuracy.
In the vicinity of such points where at least three pure imaginary eigenvalues couple, the question `Should low-order models be believed' \cite{BW09} makes sense because here a one more degree of freedom is crucial for the correct solution.

Knowledge of the rigorously established optimal solutions at every $m$ should shed light on the behavior of the optimal mass distributions and critical frequencies with the increase of the number of degrees of freedom. As a by-product, such a study will give an insight to the problem of dimension reduction and will serve a nice playground both for application of the existing methods of nonsmooth analysis and eigenvalue optimization \cite{L03, Bu06, LZ07} and for their further development in the direction of more tight relation both with the singularity theory and the applications' needs. We expect that the proposed model optimization problem will yield useful recommendations for improvement of the algorithms of optimization of real non-conservative structures in industry.

\section{Conclusions}

In the present article we proposed to consider an $m$-link Ziegler-Lobas pendulum loaded by a follower force as a model problem for studying the basic properties of optimization of circulatory systems with respect to stability criteria. The model appears to be never studied from the point of view of structural optimization although its distributed analogue --- the Beck's column --- is a popular subject for investigation in numerical optimization of non-conservative systems. We have found the optimal solutions of the two-link undamped Ziegler's pendulum and established that they correspond to a vanishing mass either at the free end or at the second joint. At these designs the critical flutter load attains its maxima. A uniform design corresponds to an absolute minimum of the flutter load. These findings are qualitatively in agreement with the properties of known numerically optimized designs of the Beck's column. For the first time we have shown that the merit functional of the two-link undamped Ziegler's pendulum depending of two masses forms a singular surface with the Whitney umbrella singularity that governs optimization.
We have formulated a problem of optimal mass distribution for an $m$-link Ziegler-Lobas pendulum and studied the case of $m=3$ degrees of freedom. It turned out that in this case there exist extrema that correspond to vanishing masses at some of the joints. However, the new degree of freedom is responsible for the new types of extrema in the vicinity of the singular points of the merit functional that correspond to triple pure imaginary eigenvalues with the Jordan block. Near such points, the stability boundary of the three-link pendulum qualitatively differs from that of a two-link one that clearly demonstrates the limits of a two-mode approximation. Finally, we have shown that although in the presence of damping the merit functional changes qualitatively, there still exist optimal solutions with some of the masses vanishing at the joints of the pendulum.

\section*{Acknowledgements}
The work has been supported by the Alexander von Humboldt Foundation.


\begin{thebibliography}{}
\bibitem{B96}
J. Bernoulli, Problema novum ad cuius solutionem Mathematici invitantur (A new problem that mathematicians are invited to solve).
Acta Eruditorum, 15 (1696) 264--269.
\end{thebibliography}
\end{document}